\input amstex 
\documentstyle{amsppt}
\magnification 1200
\NoRunningHeads
\topmatter
\title
Martingale approximation and optimality of some conditions for the central limit theorem 
\endtitle
\author
Dalibor Voln\'y
\endauthor
\affil
Universit\'e de Rouen 
\endaffil
\address
D\'epartement de Math\'ematiques, 
Universit\'e de Rouen, 
76801 Saint Etienne du Rouvray,
France
\endaddress
\email{Dalibor.Volny\@univ-rouen.fr}
\endemail
\abstract 
Let $(X_i)$ be a stationary and ergodic Markov chain with kernel $Q$, $f$ an $L^2$ function on its state space. 
If $Q$ is a normal operator and $f = (I-Q)^{1/2}g$ (which is equivalent to the convergence of
$\sum_{n=1}^\infty \frac{\sum_{k=0}^{n-1}Q^kf}{n^{3/2}}$ in $L^2$), 
we have the central limit theorem (cf\. \cite{D-L 1}, \cite{G-L 2}). Without assuming normality of $Q$, the CLT is implied by the convergence of 
$\sum_{n=1}^\infty \frac{\|\sum_{k=0}^{n-1}Q^kf\|_2}{n^{3/2}}$, in particular by $\|\sum_{k=0}^{n-1}Q^kf\|_2 =
o(\sqrt n/\log^q n)$, $q>1$ by \cite{M-Wu} and \cite{Wu-Wo} respectively. We shall show that if $Q$ is not normal and 
$f\in (I-Q)^{1/2} L^2$, or if the conditions of Maxwell and Woodroofe
or of Wu and Woodroofe are weakened to $\sum_{n=1}^\infty c_n\frac{\|\sum_{k=0}^{n-1}Q^kf\|_2}{n^{3/2}}<\infty$ for some sequence
$c_n\searrow 0$, or by $\|\sum_{k=0}^{n-1}Q^kf\|_2 = O(\sqrt n/\log n)$, the CLT need not hold.
\endabstract
\endtopmatter
\document
\subheading{1. Introduction} 
Let $(S,\Cal B,\nu)$ be a probability space, $(\xi_i)$ a homogeneous and ergodic Markov chain with state space $S$,
transition operator $Q$, and stationary distribution $\nu$. For a measurable function $g$ on $S$, $(g(\xi_i))$ is 
then a
stationary random process; we shall study the central limit theorem for 
$$
  S_n(g) = \sum_{i=0}^{n-1} g(\xi_i)
$$
where $g\in L_0^2(\nu)$, i.e\. is square integrable and has zero mean. Gordin and Lif\v sic (\cite{G-L 1}) showed that if $g$
is a solution of the equation
$$
  g = (I-Q)h = h - Qh
$$
with $h\in L^2(\nu)$ then a {\it martingale approximation} giving the CLT exists. More precisely, there exists a martingale 
difference sequence of $m(\xi_i) = h(\xi_i)-Qh(\xi_{i+1})$ such that $\|S_n(g-m)\|_2/\sqrt n \to 0$ (as shown in \cite{Vo 1}, this condition 
is equivalent to Gordin's condition from \cite{G}).
The result was extended to normal operators $Q$
and functions $g$ satisfying
$$
  g = (I-Q)^{1/2}h \tag1
$$
with $h\in L^2$. The operator $(I-Q)^{1/2}$ is defined using the series of the function $\sqrt{1-x}$, $x\in [-1,1]$
(cf\. \cite{D-L 1}). For reversible operators $Q$ the result was proved by Kipnis and Varadhan in 1986 (\cite{K-V}),
for normal operators $Q$ the result appears in 1981 in \cite{G-L 2} with a proof published later in \cite{G-L 3}, in 1996 the result was independently proved by Derriennic and Lin
in \cite{D-L 1}. Derriennic and Lin formulated the condition (1) in its present form; the other authors used 
spectral forms of the condition. As noticed by Gordin and Holzmann (\cite{G-Ho}), (1) is equivalent to the convergence of
$$
  \sum_{n=1}^\infty \frac{\sum_{k=0}^{n-1}Q^kg}{n^{3/2}}\quad\text{in}\quad L^2. \tag2
$$
They did not present a proof of this statement; for making reader's homework given in \cite{G-Ho} easier, let us give several arguments.

By \cite{D-L 1} we have that $f\in \sqrt{I-Q} L^2$ iff $\sum_{j=0}^\infty (j+1)a_{j+1} Q^jf$ converges, where $a_j = -{{j-3/2}\choose j}$;
by the Stirling formula, $a_j \sim \frac{j^{-3/2}}{2\sqrt \pi}$ with a summable error term.

From the convergence of $\sum_{k=1}^\infty \frac1{\sqrt k} Q^kf$ and Kronecker's lemma it follows \newline
$\frac1{\sqrt n} \sum_{k=1}^n Q^kf\to 0$.
\newline
By double summation and elementary estimation we can find that the series (2)
converges iff the sequence of $\sum_{k=1}^n \Big( \frac1{\sqrt k} - \frac1{\sqrt n}\Big) Q^kf$ converges.
We thus have that $f\in \sqrt{I-Q} L^2$ implies the convergence in (2).

The proof of the converse copies 
the proof of (a more general) Lemma 4.1 in \cite{Cu 2}.

Suppose that the series (2) converges. Denote $S_k(g) = \sum_{j=0}^{k-1}Q^jg$. We have
$$
  \sum_{k=n}^{2n-1} \frac{S_k(g)}{k^{3/2}} = S_n(g) \sum_{k=n}^{2n-1} \frac1{k^{3/2}} + Q^n\Big(\sum_{k=1}^{n-1} 
  \frac{S_k(g)}{(n+k)^{3/2}}\Big) \to 0.
$$
Because the sequence of partial sums of the series (2) is Cauchy and $Q$ is a Markov operator, 
$\big\|\sum_{k=1}^{n-1} \frac{S_k(g)}{(n+k)^{3/2}}\big\|_2 \to 0$ will imply that $\frac{S_n(g)}{n^{3/2}}$ converges to $0$.
Let us prove it.

Define $R_n = \sum_{k=1}^\infty \frac{S_k(g)}{k^{3/2}}$. We have
$$\gather
  \sum_{k=1}^{n-1} \frac{S_k(g)}{(n+k)^{3/2}} = \sum_{k=1}^{n-1} (R_k - R_{k+1}) \frac{k^{3/2}}{(n+k)^{3/2}} = \\
  = \sum_{k=2}^{n-1} R_k \Big(\frac{k^{3/2}}{(n+k)^{3/2}} - \frac{(k-1)^{3/2}}{(n+k-1)^{3/2}}\Big) + \frac{R_1}{(n+1)^{3/2}}
  - \frac{R_n(n-1)^{3/2}}{(2n-1)^{3/2}}. \endgather
$$
By (2), $\|R_n\|_2 \to 0$ hence given an $\epsilon >0$ there exists an $n_0$ such that $\|R_n\|_2 < \epsilon$ for 
$n\geq n_0$. We thus have
$$\multline
  \Big\|\sum_{k=2}^{n-1} R_k \Big(\frac{k^{3/2}}{(n+k)^{3/2}} - \frac{(k-1)^{3/2}}{(n+k-1)^{3/2}}\Big)\Big\|_2 \leq \\
 \leq \max_{2\leq k \leq n_0} \|R_k\|_2 \sum_{k=2}^{n_0} \Big(\frac{k^{3/2}}{(n+k)^{3/2}} - \frac{(k-1)^{3/2}}{(n+k-1)^{3/2}}\Big)
  + \epsilon, \endmultline
$$
therefore $\|\sum_{k=1}^{n-1} \frac{S_k(g)}{(n+k)^{3/2}}\big\|_2 \to 0$, hence $\frac{S_n(g)}{n^{3/2}}$. The convergence of the
series (2) is equivalent to the convergence of $\sum_{k=1}^n \Big( \frac1{\sqrt k} - \frac1{\sqrt n}\Big) Q^kg$, we thus get
that $\sum_{k=1}^n \frac1{\sqrt k} Q^kg$ converges, which is equivalent to $g\in \sqrt{I-Q} L^2$.
\bigskip

As we shall see in Theorem 1, without normality of $Q$ the condition (2) does not imply the CLT.
Maxwell and Woodroofe have shown in \cite{M-Wo} that if
$$
  \sum_{n=1}^\infty \frac{\|\sum_{k=0}^{n-1}Q^kg\|_2}{n^{3/2}} < \infty \tag3
$$
then (without any other assumptions on the Markov operator $Q$) the martingale approximation (and the CLT) takes place. 
\medskip

Let $(\Omega, \Cal A,\mu)$ be a probability space with a bijective, bimeasurable and
measure preserving transformation $T$. 
For a measurable function $f$ on $\Omega$, $(f\circ T^i)_i$ is a (strictly)
stationary process and reciprocally, using the first cannonical process, we get that any (strictly) stationary process  can be represented 
in this way:
to a stationary process $(X_i)_i$ defined on a probability space $(\Omega', \Cal A',P)$ we define a mapping $\psi$: 
$\Omega' \to \Bbb R^\Bbb Z$ by $\psi(\omega) = (X_i(\omega))_i$, and
on $\Omega=\Bbb R^\Bbb Z$ equipped with the product $\sigma$-algebra $\Cal A$ we define the image measure $\mu = P\circ \psi^{-1}$. By $T$ we
denote the left shift 
transformation of $\Omega$ onto itself, $(T\omega)_i = \omega_{i-1}$. If $Z_i$ is the projection of $\Omega =\Bbb R^\Bbb Z$ to the i-th
coordinate then the distribution of the process $(Z_i)$ is the same as the distribution of $(X_i)$ and $Z_i = Z_0\circ T^i$.

Any stationary process can be represented by a homogeneous and stationary Markov chain (\cite{Wu-Wo}, a similar idea appears already in \cite{R, p.65}). To a process $(X_i)$
we associate a Markov chain $(\xi_k)$ with $\Bbb R^\Bbb N$ for the state space, where $\xi_k = (\dots,X_{k-1},X_k)$, the transition operator $Q$ 
is given by
$Q(x,B) = \mu(\xi_1\in B|\xi_0=x) = \mu(\xi_1\in B|X_0=x_0, X_{-1}=x_{-1},\dots)$ where $x = (\dots,x_{-1},x_0)\in \Bbb R^\Bbb N$, and a stationary distribution is given by the 
distribution of the process $(X_i)_{i\leq 0}$. 
For $g(x) = x_0$, $x = (\dots,x_{-1},x_0)\in \Bbb R^\Bbb N$, the process $g(\xi_i)$ has the same distribution as $(X_i)$.

For $g$ integrable we have $Qg(\xi_i) = E(g(\xi_{i+1})|\xi_i)$. The conditions (2) and (3) thus can be expressed in the following way.

\noindent Let $(\Omega, \Cal A, \mu, T)$ be a dynamical system (a probability space with a bimeasurable and measure preserving 
bijective transformation $T : \Omega \to \Omega$), $\Cal F_i$ an increasing filtration with $T^{-1}\Cal F_i = \Cal F_{i+1}$, $f$ is a square integrable and zero mean function on $\Omega$, $\Cal F_0$-measurable. We denote
$$
  S_n(f) = \sum_{i=0}^{n-1} f\circ T^i. 
$$
The convergence in (2) is then equivalent to the convergence (in $L^2$) of
$$
  \sum_{n=1}^\infty \frac{E(S_n(f)\,|\,\Cal F_0)}{n^{3/2}} \tag2'
$$
and (3) becomes
$$
  \sum_{n=1}^\infty \frac{\|E(S_n(f)\,|\,\Cal F_0)\|_2}{n^{3/2}} < \infty. \tag3'
$$ 

\underbar{Remark 1.} Notice that in (2') and (3'), the natural filtration need not be used, it is sufficient to suppose that the process 
$(f\circ T^i)$ is adapted to $(\Cal F_i)$. The natural filtration is the smallest
filtration with respect to which the process $(f\circ T^i)$ is adapted, hence the convergence in (2'), (3') for $(\Cal F_i)$ implies
the convergence for the natural filtration.

\underbar{Remark 2.} In this article we suppose that the dynamical system $(\Omega, \Cal A,\mu, T)$ is ergodic,
i.e\. for all sets $A\in\Cal A$ such that $A=T^{-1}A$, it is $\mu(A)=0$ or $\mu(A)=1$. A stationary (here, this always means strictly
stationary) process $(X_i)$ is said to be ergodic if there exists an ergodic dynamical system with a process $(f\circ T^i)_i$ equally
distributed as $(X_i)$. Remark that an ergodic process can be represented within a non ergodic dynamical system.

Let $X_i= f\circ T^i$. Then the mapping $\psi$: $\Omega \to \Bbb R^\Bbb Z$ defined by $\psi(\omega) = (X_i(\omega))_i$ is a factor map of 
$\Omega$ onto $\Bbb R^\Bbb Z$, hence if $(\Omega, \Cal A,\mu, T)$ is ergodic then the dynamical system defined by the first cannonical
process is ergodic (cf\. \cite{C-F-S}). The process $(X_i)$  is thus ergodic if and only if the associated dynamical system defined by the first 
cannonical process is ergodic.

\underbar{Remark 3.} Let $(f\circ T^i)$ be the first cannonical process representation for a stationary process $(X_i)$ and let $(g(\xi_i))$ be the
Markov chain representation of $(f\circ T^i)_i$ ($f$ is thus the projection to the zero-th coordinate of 
$\Omega = \Bbb R^\Bbb Z$ and $T$ is the 
left shift). Ergodicity of the Markov chain $(\xi_i)$ 
is equivalent to ergodicity of $(X_i)$: \newline
We define $S=\Bbb R^\Bbb N$, 
$\Cal B$ is the product $\sigma$-algebra on $S^\Bbb Z$.
Let us define $\phi$: $\Omega =\Bbb R^\Bbb Z \to (\Bbb R^\Bbb N)^\Bbb Z=S^\Bbb Z$ by $\phi((\omega_k)_{k\in \Bbb Z}) = 
((\omega_i)_{i\leq k})_{k\in \Bbb Z}$, $\nu=\mu\circ \phi^{-1}$.
The measure $\nu=\mu\circ \phi^{-1}$ is invariant with respect to the left shift $\tau$ on $S^\Bbb Z$, $\phi$ is a bimeasurable bijection of
$\Bbb R^\Bbb N$ onto $\phi(\Bbb R^\Bbb N)\in \Cal B$ which 
commutes with the transformations $T$, $\tau$: $\phi\circ T = \tau\circ\phi$. The dynamical systems $(\Omega, \Cal A,\mu, T)$ and 
$(S^\Bbb Z,\Cal B,\nu,\tau)$ are thus isomorphic. 

\noindent If $\xi_i$: $S^\Bbb Z\to S$ are the coordinate projections, $(\xi_i)$ is a  Markov chain and for $g$: $\Bbb R^\Bbb N\to \Bbb R$, $g(x) = x_0$,
$x=(\dots,x_{-1},x_0)$, $(g(\xi_i))$ has the same distribution as $(f\circ T^i)$.

\noindent Because ergodicity is invariant with respect to isomorphism (\cite{CFS}) and the dynamical system
$(S^\Bbb Z,\Cal B,\nu,\tau)$ is the first cannonical process for
$(\xi_i)$, ergodicity of the Markov chain $(\xi_i)$ is equivalent to ergodicity of $(X_i)$.

\underbar{Remark 4.} In \cite{Vo 2} a nonadapted version of the Maxwell-Woodroofe aproximation (3') have been found. 
\medskip

In the present paper we will deal with optimality of the conditions (2') and (3'), hence also of (2) and (3). 

\proclaim{Theorem 1} There exists an ergodic process $(f\circ T^i)$ such that the series
$$
  \sum_{n=1}^\infty \frac{E(S_n(f)\,|\,\Cal F_0)}{n^{3/2}} \tag2'
$$
converges in $L^2$,
but for two different subsequences $(n_k')$, $(n_k'')$, the distributions of $S_{n_k'}/\sigma_{n_k'}$ and
$S_{n_k''}/\sigma_{n_k''}$ converge to different limits.

Therefore, if $(\xi_i)$ is a homogeneous ergodic Markov chain $(\xi_i)$ with a transition operator $Q$, then without normality 
of $Q$, the condition $g\in (I-Q)^{1/2}L^2$ is not sufficient for the CLT.
\endproclaim

In the next two theorems we show that in the central limit theorem of Maxwell and Woodroofe, the rate of convergence
of $\|E(S_n(f)|\Cal F_0)\|_2$ towards 0 is practically optimal. We denote $\sigma_n = \|ES_n(f)\|_2$.

\proclaim{Theorem 2} For any sequence of positive reals $c_n\to 0$ there exists an ergodic process $(f\circ T^i)$ such that 
$$
  \sum_{n=1}^\infty c_n\frac{\|E(S_n(f)|\Cal F_0)\|_2}{n^{3/2}} < \infty \tag4
$$
but for two different subsequences $(n_k')$, $(n_k'')$, the distributions of $S_{n_k'}/\sigma_{n_k'}$ and
$S_{n_k''}/\sigma_{n_k''}$ converge to different limits.
\endproclaim

\proclaim{Theorem 3} There exists an ergodic process $(f\circ T^i)$ such that 
$$
  \|E(S_n(f)|\Cal F_0)\|_2 = o\Big(\frac{\sqrt n}{\log n}\Big),  \tag5 
$$
but for two different subsequences $(n_k')$, $(n_k'')$, the distributions of $S_{n_k'}/\sigma_{n_k'}$ and
$S_{n_k''}/\sigma_{n_k''}$ converge to different limits.
\endproclaim

\underbar{Remark 5.} In \cite{Pe-U}, under the assumption (4) Peligrad and Utev have shown that
there exists an $f$
such the sequence of $S_n(f)/\sqrt n$ is not stochastically bounded. 

\underbar{Remark 6.} In \cite{Wu-Wo} M. Woodroofe and W.B. Wu showed that if $\|E(S_n(f)|\Cal F_0)\|_2 = o(\sigma_n)$ then $\sigma_n = h(n)\sqrt n$ where $h(n)$ is a slowly varying function in the sense of Karamata and there exists
an array $D_{n,i}$ of martingale differences such that for each $n$, the sequence $(D_{n,i})_i$ is a strictly 
stationary martingale difference sequence and $\|S_n(f) - \sum_{i=0}^{n-1} D_{n,i}\|_2 = o(\sigma_n)$. In particular,
under the condition (4) this approximation still takes place. 

In the same paper \cite{Wu-Wo}, Wu and Woodroofe used the condition
$$
  \|E(S_n(f)\,|\,\Cal F_0)\|_2 = o\Big(\frac{\sqrt n}{\log^qn}\Big) \tag6
$$
where $q>1$. (6) implies (3') hence a martingale approximation and the CLT. Theorem 3 shows that the result cannot be extended to $q=1$.

\underbar{Remark 7.} The condition (5) which can be written as $\|\sum_{k=0}^{n-1}Q^kg\|_2 = o\big(\frac{\sqrt n}{\log n}\big)$ 
implies that $\sum_{n=1}^\infty \frac{\|\sum_{k=0}^{n-1}Q^kg\|_2^2}{n^{2}} < \infty$. 
By \cite{Cu 1}, Proposition 2.2, for a normal operator $Q$ this implies $g\in \sqrt{I-Q}L^2$ hence a CLT; 
Theorem 3 shows that the assumption of normality cannot be lifted.

\underbar{Remark 8.} From the construction it follows that in Theorems 1-3, the variances $\sigma_n^2$ of $S_n(f)$ grow faster than linearly.
It thus remains an open problem whether with a supplementary assumption $\sigma_n^2/n\to const.$ the CLT would hold.
As shown by a counter example in \cite{Kl-Vo2}, this assumption is not sufficient for
$q\leq 1/2$ (in \cite{Kl-Vo2}, a function $f$ is found such that (6) holds with $q= 1/2$, $\sigma_n^2/n\to const.$, but the distributions of $S_n(f)/\sqrt n$ do not converge); 
the only exponents to consider are thus $1/2<q\leq 1$. 

It also remains an open question whether the CLT would hold for $f\in L^{2+\delta}$ for some $\delta>0$.

\subheading{2. Proofs}
We first define an ergodic dynamical system $(\Omega,\Cal A,\mu,T)$ where the processes $(f\circ T^i)$ which we need can be found.

For $l=1,2,\dots$ we define $A_l = \{-1,0,1\}$ if $l$ is odd and $A_l=\Bbb R$ if $l$ is even, the sets $A_l$ are equipped with Borel
$\sigma$-algebras and probability measures
$\nu_l$ such that $\nu_l(\{-1\}) = 1/(2N_l) = \nu_l(\{1\})$ and $\nu_l(\{0\}) = 1-1/N_l$ for $l$ odd, $\nu_l = \Cal N(0,1)$ for $l$ even.
For each $i\in\Bbb Z$ we define $\Omega_i = \underset l=1 \to{\overset \infty \to{\times}} A_l$; $\Omega_i$ is equipped with the 
product measure $\mu_i= \underset l=1 \to{\overset \infty \to{\otimes}} \nu_l$.
On the set $\Omega = \underset i\in\Bbb Z \to{\times} \Omega_i$, $\mu$ is
the product measure $\mu = \underset i\in\Bbb Z \to{\otimes} \mu_i$ and $T$ is the left shift
transformation. By definition, the dynamical system $(\Omega,\Cal A,\mu,T)$ is Bernoulli hence ergodic.
By $\Cal F_k$, $k\in \Bbb Z$, we denote the $\sigma$-algebra generated by projections of $\Omega$ onto $\Omega_l$, $l\leq k$.

In all of the text, $\log$ will denote the dyadic logarithm. By $U$ we denote the unitary operator on $L^2$ defined
by $Uf = f\circ T$. $P_i$ denotes the orthogonal projection on the Hilbert space $L^2(\Cal F_i)\ominus 
L^2(\Cal F_{i-1})$, i.e\. 
$$
  P_if = E(f|\Cal F_i) - E(f|\Cal F_{i-1}),
$$
$f\in L^2$, $i\in \Bbb Z$.

For $k=1,2,\dots$ let $n_k = 2^k$, 
$$
  0\leq a_k\leq 1,\quad  a_k/a_{k+1} \to 1,\quad a_{k+1}\leq a_k, \quad \sum_{k=1}^\infty a_k/k =\infty;
$$
$(N_l)$, $l=1,2,\dots$, is an increasing sequence of positive integers such that 
$$
  2^{2^{l}}-1 < \sum_{k\geq 1: N_{l-1}< n_k\leq N_l} \frac{a_k}{k} < 2^{2^{l}}+1.
$$
Denote by $\pi_0$ the projection of $\Omega$ onto $\Omega_0$ and by $p_l$ the projection of $\Omega_0$ onto $A_l$. For 
$N_{l-1}< n_k\leq N_l$ we then define $e_k = \frac{a_k}{k}p_l\circ \pi_0$. Then
$$
  \|e_k\|_2 = a_k/k,
$$ 
and for $i\neq j$, 
$U^ie_k$ and $U^je_l$ are independent.
Notice that for $N_{l-1}< n_k\leq N_l$ the random variables $e_k$ are multiples one of
another and are independent of any $e_j$ with $n_j\leq N_{l-1}$ or $n_j > N_l$. 
In general, $e_{k'}$, $e_{k''}$ are not orthogonal but 
for all $1\leq k',k''$ it is $E(e_{k'}e_{k''})\geq 0$ and
$$
  \Big\| \sum_{k=1}^n e_k \Big\|_2 \nearrow \infty \quad\text{as}\quad n\to \infty;
$$
for $b'_k\geq b_k\geq 0$ we thus have 
$$
  \Big\| \sum_{k=1}^n b'_ke_k \Big\|_2 \geq \Big\| \sum_{k=1}^n b_ke_k \Big\|_2. \tag7
$$
\bigskip

Let
$$
  f = \sum_{k=1}^\infty \frac1{n_k} \sum_{i=0}^{n_k-1} U^{-i}e_k.
$$
We have $\|f\|_2 \leq \sum_{k=1}^\infty \|e_k\|_2/\sqrt{n_k} <\infty$ due to the exponential growth of the $n_k$s.

For a positive integer $N$ we have
$$
  S_N(f) = \sum_{k=1}^\infty \sum_{j=0}^{N-1} \sum_{i=0}^{n_k-1} \frac1{n_k} U^{j-i}e_k = 
  S'_N(f) + S''_N(f)
$$
where
$$
  S'_N(f) = S_N(f) - E(S_N(f)|\Cal F_0) = \sum_{k=1}^\infty \sum_{j=0}^{N-1} 
  \sum_{i=0}^{(j\wedge n_k)-1} 
  \frac1{n_k} U^{j-i}e_k
$$
($j\wedge n_k = \min\{j, n_k\}$) and
$$
  S''_N(f) = E(S_N(f)|\Cal F_0) = \sum_{k=1}^\infty \sum_{j=0}^{N-1} \sum_{i=j}^{n_k-1} 
  \frac1{n_k} U^{j-i}e_k.
$$
We will study the asymptotic behaviour of $S''_N(f) = E(S_N(f)|\Cal F_0)$ and $S'_N(f) = S_N(f) - E(S_N(f)|\Cal F_0)$
separately. In Lemmas 1-5 we find estimates for $E(S_N(f)|\Cal F_0)$ and an approximation of $S'_N(f)$ by sums of martingale differences.
These lemmas do not depend on distributions of the random variables $e_k$. Then we use the distributions of the random variables $e_k$ to get
a limit behaviour of the distributions of $S_n(f)/\sigma_n$ which we need.

\bigskip
\centerline{\it 1. Asymptotics of $S''_N(f) = E(S_N(f)|\Cal F_0)$.}
\medskip

\proclaim{Lemma 1}  The series
$$
  \sum_{n=1}^\infty \frac{E(S_n(f)\,|\,\Cal F_0)}{n^{3/2}} \tag2'
$$
converges in $L^2$.
\endproclaim

\demo{Proof} For $N\leq n_k$ we have
$$
  \sum_{j=0}^{N-1} \sum_{i=j}^{n_k-1} U^{j-i}e_k = \sum_{i=0}^{n_k-N} N U^{-i}e_k + 
  \sum_{i=n_k-N+1}^{n_k-1}
  (n_k-i) U^{-i}e_k
$$
and for $N>n_k$ we have
$$
  \sum_{j=0}^{N-1} \sum_{i=j}^{n_k-1} U^{j-i}e_k = \sum_{j=0}^{n_k-1} (n_k-j) U^{-j}e_k,
$$
hence
$$\multline
  S''_N(f) = \sum_{k\geq 1: n_k< N} \sum_{j=0}^{n_k-1} \frac{n_k-j}{n_k} U^{-j}e_k + \\
  + \sum_{k\geq 1: n_k\geq N} \Big[\sum_{j=0}^{n_k-N} \frac{N}{n_k} U^{-j}e_k +
  \sum_{j=n_k-N+1}^{n_k-1} \frac{n_k-j}{n_k} U^{-j}e_k \Big].
  \endmultline \tag8
$$
To prove the lemma it is thus sufficient to show that the sums
$$
  \sum_{N=1}^\infty \sum_{k\geq 1: n_k< N} \sum_{j=0}^{n_k-1} 
  \frac1{N^{3/2}}\frac{n_k-j}{n_k} U^{-j}e_k, \tag{8a}
$$

$$
  \sum_{N=1}^\infty \sum_{k\geq 1: n_k\geq N} \frac1{N^{3/2}} \sum_{j=0}^{n_k-N} 
  \frac{N}{n_k} U^{-j}e_k,  \tag{8b}
$$
$$
  \sum_{N=1}^\infty \sum_{k\geq 1: n_k\geq N} \frac1{N^{3/2}}\sum_{j=n_k-N+1}^{n_k-1} \frac{n_k-j}{n_k} U^{-j}e_k \tag{8c}
$$
converge in $L^2$.

Recall that $n_k=2^k$. We show that the sequence of partial sums for the first series is Cauchy.
Let $1\leq p<q<\infty$.

$$\gather
  \Big\|\sum_{N=p}^q \sum_{k\geq 1: 2^k< N} \sum_{j=0}^{2^k-1} 
  \frac1{N^{3/2}}\frac{2^k-j}{2^k} U^{-j}e_k \Big\|_2^2 = \\
  \Big\|\sum_{k=1}^\infty \Big(\sum_{N=(2^k+1)\vee p}^q \frac1{N^{3/2}}\Big)
  \sum_{j=0}^{2^k-1} \frac{2^k-j}{2^k} U^{-j}e_k \Big\|_2^2 = \\
  \Big\| \sum_{k=1}^\infty b_k(p,q) \sum_{j=0}^{2^k-1} \frac{2^k-j}{2^{3k/2}} U^{-j}e_k \Big\|_2^2 
  =
  \Big\| \sum_{j=0}^\infty \sum_{k\geq 1\vee\log (j+1)} b_k(p,q) \frac{2^k-j}{2^{3k/2}} U^{-j}e_k \Big\|_2^2
  \endgather
$$
where $b_k(p,q) = 2^{k/2}\sum_{N=(2^k+1)\vee p}^q \frac1{N^{3/2}}$. Notice that $b_k(p,q)$ are uniformly bounded and 
for each $k$, $\sup_{q\geq p} b_k(p,q) \searrow 0$ as $p\to \infty$. Denote $B = \sup_{k,p,q} b_k(p,q) < \infty$.
Using $\|e_k\|_2\leq 1/k$ we deduce

$$\gather
  \Big\| \sum_{j=0}^\infty \sum_{k\geq 1\vee\log (j+1)} b_k(p,q) \frac{2^k-j}{2^{3k/2}} U^{-j}e_k \Big\|_2^2
  \leq 
  \sum_{j=0}^\infty \Big\| \sum_{k\geq 1\vee\log (j+1)} \frac{b_k(p,q)}{2^{k/2}} U^{-j}e_k \Big\|_2^2 \leq \\
  \leq  \sum_{j=0}^\infty \Big(\sum_{k\geq 1\vee\log (j+1)} \frac{b_k(p,q)}{k2^{k/2}}  \Big)^2
  \endgather
$$
where
$$
  \sum_{k\geq 1\vee\log (j+1)} \frac{b_k(p,q)}{k2^{k/2}}  \leq
  \frac{\sqrt 2}{\sqrt 2 -1} \frac{B}{\sqrt{(j+1)}[1\vee\log (j+1)]},
$$
hence 
$$
  \sum_{j=0}^\infty \Big(\sum_{k\geq 1\vee\log (j+1)} \frac{b_k(p,q)}{k2^{k/2}}  \Big)^2 < \infty,
$$
therefore
$$\multline
  \Big\| \sum_{j=0}^\infty \sum_{k\geq 1\vee\log (j+1)} b_k(p,q) \frac{2^k-j}{2^{3k/2}} U^{-j}e_k \Big\|_2^2
  \leq \\
  \leq \lim_{p\to\infty}\sup_{q\geq p} \sum_{j=0}^\infty \Big(\sum_{k\geq 1\vee\log (j+1)} \frac{b_k(p,q)}{k2^{k/2}}  \Big)^2 = 0.
  \endmultline
$$

For the second sum we define 
$$
  b_k(p,q) = \frac1{2^{k/2}} \sum_{N=p}^{2^k\wedge q} \frac{1}{N^{1/2}}
$$
and get
$$\gather
  \Big\|\sum_{N=p}^q \sum_{k\geq 1: 2^k\geq N} \frac1{N^{3/2}} \sum_{j=0}^{2^k-N} 
  \frac{N}{2^k} U^{-j}e_k \Big\|_2^2 = 
  \Big\|\sum_{N=p}^q \sum_{j=0}^\infty \sum_{k\geq \log (N+j)} \frac{1}{N^{1/2}2^k} 
  U^{-j}e_k \Big\|_2^2 = \\
  = \Big\| \sum_{j=0}^\infty \sum_{k\geq \log(j+1)} \Big(\sum_{N=p}^{q\wedge (2^k-j)} 
  \frac{1}{N^{1/2}}\Big) \frac{1}{2^k} U^{-j}e_k \Big\|_2^2 \leq \\
  \leq \sum_{j=0}^\infty \Big\| \sum_{k\geq \log(j+1)} \frac{b_k(p,q)}{2^{k/2}} U^{-j}e_k \Big\|_2^2 \to 0
\endgather
$$ 
as $p\to \infty$ using similar arguments as in the preceding case.

For the third sum we get
$$\gather
  \Big\| \sum_{N=p}^q \sum_{k\geq 1: 2^k\geq N} \frac1{N^{3/2}} 
  \sum_{j=2^k-N+1}^{2^k-1} \frac{2^k-j}{2^k} U^{-j}e_k \Big\|_2^2 = \\
  \sum_{j=0}^\infty \Big\|\sum_{k\geq \log(j+1)} 
  \sum_{N=p\vee (2^k-j+1)}^{q\wedge 2^k}\frac1{N^{3/2}}\frac{2^k-j}{2^k} U^{-j}e_k \Big\|_2^2 \leq \\
  \sum_{j=0}^\infty \Big\|\sum_{k\geq \log(j+1)} 
  \Big(\sum_{N=p\vee (2^k-j+1)}^{q\wedge 2^k}\frac{1}{N^{1/2}}\Big) \frac{1}{2^k} U^{-j}e_k \Big\|_2^2 \to 0
\endgather
$$
as $p\to \infty$.
The sequences of partial sums for (8a), (8b), (8c) are Cauchy hence converge in $L^2$. Therefore, the series 
$\sum_{N=1}^\infty \frac{E(S_N(f)\,|\,\Cal F_0)}{N^{3/2}}$ converges in $L^2$.
\enddemo \qed

\proclaim{Lemma 2} There exists a constant $0<c<\infty$ such that 
$$\|E(S_N(f)|\Cal F_0)\|_2^2 \leq c
  \frac{Na^2_{[\log N]}}{\log^2 N}\tag9
$$
for all $N\geq 2$.
\endproclaim

\demo{Proof} We have
$$\gathered
  (1/\sqrt 6)\|e_k\|_2\sqrt{n_k}\leq \Big\|\sum_{j=0}^{n_k-1} \frac{n_k-j}{n_k} 
  U^{-j}e_k\Big\|_2 \leq \|e_k\|_2\sqrt{n_k},
  \quad n_k<N,\\
  \Big\|\sum_{j=0}^{n_k-N} \frac{N}{n_k} U^{-j}e_k + \sum_{j=n_k-N+1}^{n_k-1} 
  \frac{n_k-j}{n_k} U^{-j}e_k \Big\|_2
  \leq \frac{N}{\sqrt{n_k}} \|e_k\|_2, \quad n_k\geq N.
  \endgathered \tag10
$$

Recall that by $[x]$ we denote the integer part of $x$.
Because $n_k = 2^k$ grows exponentially fast, the norms $\|e_k\|_2$ are decreasing, and $\|e_k\|_2/\|e_{k+1}\|_2 \to 1$, there exists 
a constant $0<c<\infty$ not 
depending on $N$ such that 
$$\gather
  \sum_{k\geq 1: n_k\geq N} (N/\sqrt{n_k})\|e_k\|_2 \leq 
  c\sqrt N \|e_{[\log N]}\|_2, \\
  \sum_{k\geq 1: n_k< N} \|e_k\|_2\sqrt{n_k} \leq c \sqrt N
  \|e_{[\log N]}\|_2.
  \endgather
$$ 
Using (8) and (10) we deduce that for some constant $c>0$ we have
$$
  \|E(S_N(f)|\Cal F_0)\|_2^2 \leq cN \|e_{[\log N]}\|_2^2.
$$
Because $\|e_k\|_2 = a_k/k$,
$$
  \|E(S_N(f)|\Cal F_0)\|_2^2 \leq c \frac{Na^2_{[\log N]}}{\log^2 N}.
$$
\enddemo \qed

Recall
$$
  \|e_k\|_2 = a_k/k,\quad 0\leq a_k\leq 1,\quad \sum_{k=1}^\infty a_k/k =\infty.
$$ 
As a coroallary to Lemma 2 we get

\proclaim{Lemma 3}
$$
  \|E(S_n(f)|\Cal F_0)\|_2 = O\Big(\frac{\sqrt n}{\log n}\Big),
$$
if $a_n\to 0$ then
$$
  \|E(S_n(f)|\Cal F_0)\|_2 = o\Big(\frac{\sqrt n}{\log n}\Big).
$$

\endproclaim


\proclaim{Lemma 4} Let $c_n$ be positive real numbers, $c_n\searrow 0$. If
$$
  \sum_{n=1}^\infty \frac{a_nc_n}{n} <\infty\quad \text{then}\quad
  \sum_{n=1}^\infty c_n\frac{\|E(S_n(f)|\Cal F_0)\|_2}{n^{3/2}} <\infty.
$$
\endproclaim

\demo{Proof}
Recall Lemma 2 and $a_{k+1}\leq a_k$, denote the constant in (9) by $C$. We have
$$\multline
  \sum_{n=2}^\infty c_n\frac{\|E(S_n(f)|\Cal F_0)\|_2}{n^{3/2}} \leq C\sum_{n=2}^\infty c_n\frac{a_{[\log n]+1}}{n[\log n]}
  \leq C \sum_{k=1}^\infty \sum_{i=0}^{2^k-1} c_{2^k+i} \frac{a_k}{k(2^k+i)} \leq \\
  \leq C \sum_{k=1}^\infty c_{2^k} \frac{a_k}{k} \leq C \sum_{k=1}^\infty c_{k} \frac{a_k}{k} <\infty.
  \endmultline
$$
\enddemo \qed

\bigskip
\centerline{\it 2. Asymptotics of $S'_N(f) = S_N(f) - E(S_N(f)|\Cal F_0)$.}
\medskip

Denote $b(N) = \|\sum_{k\geq 1: n_k\leq N} e_k\|_2 = \|\sum_{k=1}^{[\log N]} e_k\|_2$; notice that by the assumptions on $e_k$, 
$b(N)\nearrow \infty$.

\proclaim{Lemma 5} 
We have
$$
  \lim_{N\to\infty} \frac1{b(N)\sqrt N} \Big\|S'_N(f) - \sum_{l=0}^{N-1} U^l \sum_{k\geq 1: n_k\leq N} e_k\Big\|_2 = 0. \tag{11}
$$
\endproclaim

\demo{Proof} Recall that
$$
  S'_N(f) =  \sum_{k=1}^\infty \sum_{j=0}^{N-1}  \sum_{i=0}^{(j\wedge n_k)-1} \frac1{n_k} U^{j-i}e_k.
$$
For $N\leq n_k$ we have
$$
  \sum_{j=0}^{N-1} \sum_{i=0}^{(j\wedge n_k)-1} U^{j-i}e_k = \sum_{j=1}^{N-1} (N-j)U^je_k 
  \tag{12}
$$
and for $N>n_k$ we have
$$
  \sum_{j=0}^{N-1} \sum_{i=0}^{(j\wedge n_k)-1} U^{j-i}e_k = 
  \sum_{j=1}^{N-n_k} n_kU^je_k + \sum_{j=N-n_k+1}^{N-1} (N-j)U^je_k. \tag{13}
$$
For all $k\geq 1$ we have $P_lU^je_k = 0$ if $j\neq l$, $P_lU^le_k = U^le_k$.
For $l\geq N$ and $l\leq 0$ we thus have $P_lS_N(f) = 0$ and for $1\leq l\leq N-1$ we, using (12) and (13), deduce
$$
  P_lS_N(f) = \sum_{k\geq 1: n_k\leq N-l} U^le_k + \sum_{k\geq 1: n_k\geq N+1-l} \frac{N-l}{n_k} U^le_k. \tag{14}
$$
Recall that $[x]$ denotes the integer part of $x$. We have
$$
  S'_N(f) = \sum_{l=1}^{N-1} P_lS_N(f) = \sum_{l=1}^{[N(1-\epsilon)]} P_lS_N(f) +
  \sum_{l=[N(1-\epsilon)]+1}^{N} P_lS_N(f)
$$
and using (14) we get
$$\multline
  \sum_{l=1}^{[N(1-\epsilon)]} P_lS_N(f) = \\
  \sum_{l=1}^{[N(1-\epsilon)]} U^l \sum_{k\geq 1: n_k\leq N-l} e_k +
  \sum_{l=1}^{[N(1-\epsilon)]} U^l \sum_{k\geq 1: n_k\geq N+1-l} \frac{N-l}{n_k} e_k = \\
  = \sum_{l=1}^{[N(1-\epsilon)]} U^l \sum_{k\geq 1: n_k\leq \epsilon N} e_k +
  \sum_{l=1}^{[N(1-\epsilon)]} U^l \sum_{k\geq 1: \epsilon N<n_k\leq N-l} e_k + \\
  + \sum_{l=1}^{[N(1-\epsilon)]} U^l \sum_{k\geq 1: n_k\geq N+1-l} \frac{N-l}{n_k} e_k.
  \endmultline
$$
Because $n_k=2^k$, 
$$
  \Big\|\sum_{k\geq 1: n_k\geq N+1-l} \frac{N-l}{n_k} U^le_k\Big\|_2 \leq 2\|e_{[\log (N-l)]}\|_2 \leq 2/\log (N-l), \tag{15}
$$ 
hence
$$
  \Big\|\sum_{l=1}^{[N(1-\epsilon)]} U^l \sum_{k\geq 1: n_k\geq N+1-l} \frac{N-l}{n_k} e_k\Big\|_2 \leq \frac{\sqrt N}{\log N +\log\epsilon}
$$
(we can suppose $N$ big enough to have $\log N >|\log\epsilon|$);
$\epsilon N <n_k\leq N$ if and only if $\log N + \log\epsilon <k\leq \log N$. We thus deduce that for $\epsilon>0$ fixed and
$b(N) = \|\sum_{k=1}^{\log N} e_k\|_2\nearrow \infty$ ,
$$
  \lim_{N\to\infty} \frac1{b(N)} \Big\|\sum_{\epsilon N<n_k\leq N} e_k\Big\|_2 = 0
$$
and
$$
  \lim_{N\to\infty} \frac1{b(N)\sqrt N} \Big\|\sum_{l=1}^{[N(1-\epsilon)]} P_lS_N(f) - \sum_{l=1}^{[N(1-\epsilon)]} U^l
  \sum_{k\geq 1: n_k\leq \epsilon N} e_k\Big\|_2 = 0.
$$
For all $[N(1-\epsilon)]+1\leq l\leq N-1$ we have, by (14) and (15),
$\|P_lS_N(f)\|_2 \leq b(N-l) + 2/\log (N-l)$ hence 
$$
  \lim_{\epsilon\searrow 0}\lim_{N\to\infty} \frac1{b(N)\sqrt N} \Big\|\sum_{l=[N(1-\epsilon)]+1}^{N-1} P_lS_N(f) - 
  \sum_{l=[N(1-\epsilon)]+1}^{N-1} \sum_{k\geq 1: n_k\leq \epsilon N} U^le_k\Big\|_2 = 0
$$
and (11) follows.
\enddemo \qed

Recall that 
$$
  \|e_k\|_2 = a_k/k,\quad 0\leq a_k\leq 1,\quad \sum_{k=1}^\infty a_k/k =\infty
$$
$N_l$, $l=1,2,\dots$, is an increasing sequence of positive integers such that 
$$
  2^{2^{l}}-1 < \sum_{k\geq 1: N_{l-1}< n_k\leq N_l} \frac{a_k}{k} < 2^{2^{l}}+1; \tag{16}
$$
for $N_{l-1}< n_k\leq N_l$ the random variables $e_k$ are multiples one of
another and are independent of any $e_j$ with $n_j\leq N_{l-1}$ or $n_j > N_l$. 

\proclaim{Lemma 6}
Along $l$ odd the distributions of
$$
  \frac1{b(N_l)\sqrt N_l} \sum_{j=0}^{N_l-1} U^j \Big(\sum_{k\geq 1: N_{l-1}<n_k\leq N_l} e_k\Big)
$$
weakly converge to the symmetrised Poisson distribution with parameter $\lambda=1/2$ and for $l$ even to the standard normal
distribution.
\endproclaim

\demo{Proof}
From the definition of the functions $e_k$ it follows that $b^2(N_l) = b^2(N_{l-1}) + \|\sum_{k\geq 1: N_{l-1}<n_k\leq N_l} e_k\|_2^2$ 
and $\|\sum_{k\geq 1: N_{l-1}<n_k\leq N_l} e_k\|_2
= \sum_{k\geq 1: N_{l-1}<n_k\leq N_l} a_k/k \sim 2^{2^{l}}$ hence $\|\sum_{k\geq 1: N_{l-1}<n_k\leq N_l} e_k\|_2 \sim b(N_l)$.

For $l$ odd, the random variable $\sum_{k\geq 1: N_{l-1}<n_k\leq N_l} e_k$ takes values $\approx \pm b(N_l)\sqrt N_l$
with probabilities $1/(2N_l)$ and 0 with probability $1-1/N_l$.

For $l$ even, $\sum_{k\geq 1: N_{l-1}<n_k\leq N_l} e_k$ is normally distributed with zero mean and variance $\approx b^2(N_l)$.



By the assumptions, $U^i \sum_{k\geq 1: N_{l-1}<n_k\leq N_l} e_k$ are independent random variables and the statement of the lemma follows.
\enddemo \qed

\demo{Proof of Theorem 1 and Theorem 3}
By Lemma 1 the series $\sum_{n=1}^\infty \frac{E(S_n(f)\,|\,\Cal F_0)}{n^{3/2}}$ converges in $L^2$.
We can have $a_n\searrow 0$, e.g\. $a_n =1/\log n$ for $n\geq 2$; by Lemma 3 then $\|E(S_n(f)|\Cal F_0)\|_2 = 
o\Big(\frac{\sqrt n}{\log n}\Big)$. From Lemma 3 and Lemma 5 it follows
$$
  \Big\|E(S_n(f) - \sum_{l=0}^{N-1} U^l \sum_{k\geq 1: n_k\leq N} e_k\Big\|_2 = o\Big(\frac1{b(N)\sqrt N}\Big),
$$
using Lemma 6 we thus get that along $l$ odd the distributions of $\frac1{b(N)\sqrt N}S_n(f)$ weakly converge to the symmetrised Poisson 
distribution with parameter $\lambda=1/2$ and for $l$ even to the standard normal law.
\enddemo \qed

\demo{Proof of Theorem 2}
Let $c_k>0$; without loss of generality we can suppose $c_k\searrow 0$. We define $k_0=1$ and for $n=1,2,\dots$ let $k_n$ be 
the first $k$ such that $c_k\leq 1/2^n$ and $k_n-k_{n-1} \geq n$; $a_{k_0}=1$ and for $n\geq 1$, $a_{k_n}$ is the minimum of
$a_{k_{n-1}}$ and $(\sum_{j= k_{n-1}+1}^{k_n} \frac{1}{j})^{-1}$. For $k_{n-1}+1\leq j\leq k_n-1$ we define 
$$
  a_j = a_{k_{n-1}} + \alpha_j(a_{k_{n}}-a_{k_{n-1}}), \quad \alpha_j = \cases \frac{j-k_{n-1}}n \,\,\,\,\text{if}\,\,\,\, 
                                                                                          1\leq j-k_{n-1} \leq n \\
                                         1 \,\,\,\,\text{if}\,\,\,\, n+1 \leq j-k_{n-1}\leq k_n-k_{n-1} \endcases
$$
Then 
$$
  \sum_{n=1}^\infty \frac{a_nc_n}{n} <\infty, \quad \sum_{n=1}^\infty \frac{a_n}{n} = \infty.
$$

To verify the first inequality we notice that
$$
  \sum_{n=1}^\infty \frac{a_nc_n}{n} \leq \sum_{n=1}^\infty \frac1{2^n} \Big[a_{k_n} \sum_{j=k_{n-1}+1}^{k_n} \frac1j +
  (a_{k_{n-1}}-a_{k_n}) \sum_{j=k_{n-1}+1}^{n} \frac{1-\alpha_j}{j} \Big].
$$
By definition, $a_{k_n} \sum_{j=k_{n-1}+1}^{k_n} \frac1j \leq 1$, and by boundedness of $a_n$ we have $\sum_{n=1}^\infty 
\frac{na_{k_{n-1}}}{2^n} <\infty$.


To verify the second inequality we notice
$$
  \sum_{n=1}^\infty \frac{a_n}{n} \geq \sum_{n=1}^\infty a_{k_n} \sum_{j=k_{n-1}+1}^{k_n} \frac1j.
$$
If $a_{k_n} = (\sum_{j= k_{n-1}+1}^{k_n} \frac{1}{j})^{-1}$ for infinitely many $n$ then the sum is infinite.
Otherwise, from some $n_0$ on, the sequence $a_n$ is constant and strictly positive; the series is infinite as well.

A simple calculation shows that the sequence of $a_n$ is decreasing, $0\leq a_n\leq 1$, and $a_n/a_{n+1} \to 1$.

The inequality $\sum_{n=1}^\infty \frac{a_nc_n}{n} <\infty$ makes the assumptions of Theorem 2 satisfied and the equality $\sum_{n=1}^\infty 
\frac{a_n}{n} = \infty$, together with Lemma 6, implies that $f$ can be defined so that
along $l$ odd the distributions of $\frac1{b(N)\sqrt N}S_n(f)$ weakly converge to the symmetrised Poisson 
distribution with parameter $\lambda=2$ and for $l$ even to the standard normal law.
\enddemo \qed

\bigskip
{\bf Acknowledgement.} The author thanks the referee for having pointed out several mistakes in preceding versions of the article
and for many helpful comments and remarks. I thank Dr\. Christophe Cuny for having let me use a proof from his (yet) unpublished
paper \cite{Cu 2}.


\Refs
\widestnumber\key{Wu-Wo}
\ref \key B \by Billingsley, P. \paper The Lindeberg-L\'evy theorem for martingales \jour Proc. Amer. Math. Soc.
\vol 12 \pages 788-792 \yr 1961 \endref 
\ref \key C-F-S \by Cornfeld, I.P., Fomin, S.V., and Sinai, Ya.G. \book Ergodic Theory \publ Springer \publaddr Berlin \yr 1982 \endref
\ref \key Cu 1 \by Cuny, Ch. \paper Pointwise ergodic theorems with rate and application to limit theorems for stationary processes
\paperinfo sumitted for publication, arXiv:0904.0185 \yr 2009 \endref
\ref \key Cu 2 \by Cuny, Ch. \paper Norm convergence of some power-series of operators in $L^p$ with applications in ergodic theory
\paperinfo sumitted for publication \yr 2009 \endref
\ref \key D-L 1 \by Derriennic, Y. and Lin, M. \paper Sur le th\'eor\`eme limite central de Kipnis et Varadhan pour les cha\^\i nes 
r\'eversibles ou normales \jour CRAS 323 \pages 1053-1057 \yr 1996 \endref 
\ref \key D-L 2 \by Derriennic, Y and Lin, M. \paper The central limit theorem for Markov
chains with normal transition operators, started at a point \jour Probab. Theory Relat. Fields
\vol 119 \pages 509-528 \yr 2001 \endref
\ref \key G \by Gordin, M.I. \paper A central limit theorem for stationary processes \jour Soviet Math. Dokl. \vol 10 
\pages 1174-1176 \yr 1969 \endref
\ref \key G-Ho \by Gordin, M.I. and Holzmann, H. \paper The central limit theorem for stationary Markov chains under invariant splittings 
\jour Stochastics and Dynamics \vol 4 \pages 15-30 \yr 2004 \endref 
\ref \key G-L 1 \by Gordin, M.I. and Lif\v sic, B.A. \paper Central limit theorem for statioanry processes \jour Soviet Math. Doklady 19 
\pages 392-394 \yr 1978 \endref 
\ref \key G-L 2 \by Gordin, M.I. and Lif\v sic, B.A. \paper A remark about a Markov process with normal transition operator
\paperinfo In: Third Vilnius Conference on Probability and Statistics \vol 1 \pages 147-148 \yr 1981 \endref
\ref \key G-L 3 \by Gordin, M.I. and Lif\v sic, B.A. \paper The central limit theorem for Markov processes with normal transition operator,
and a strong form of the central limit theorem \paperinfo §IV.7 and §IV.8 in Borodin and Ibragimov, Limit theorems for functionals of
random walks, Proc. Steklov Inst. Math. 195(1994), English translation AMS (1995) \endref
\ref \key Ha-He \by Hall, P. and Heyde, C.C. \book Martingale Limit Theory           
and its Application \publ Academic Press \publaddr New York \yr 1980 \endref 
\ref \key I \by Ibragimov, I.A.  \paper A central limit theorem for a class of dependent random variables 
\jour Theory Probab. Appl. \vol 8 \pages 83-89 \yr 1963 \endref
\ref \key K-V \by Kipnis, C. and Varadhan, S.R.S. \paper Central limit theorem for additive functionals of reversible
Markov processes and applications to simple exclusions \jour Comm. Math. Phys. \vol 104 \pages 1-19 \yr 1986 \endref
\ref \key Kl-Vo 1 \by Klicnarov\'a, J. and Voln\'y, D. \paper An invariance principle for non adapted processes
\jour C.R. Acad. Sci. Paris Ser 1 \vol 345/5 \pages 283-287 \yr 2007 \endref
\ref \key Kl-Vo 2 \by Klicnarov\'a, J. and Voln\'y, D. \paper Exactness of a Wu-Woodroofe's approximation with linear growth of variances \jour Stoch. Proc. and their Appl. \vol 119 \pages 2158-2165 \yr 2009 \endref
\ref \key M-Wo \by Maxwell, M. and Woodroofe, M. \paper Central limit theorems for additive
functionals of Markov chains \jour Ann. Probab. \vol 28 \pages 713-724 \yr 2000 \endref
\ref \key P-U \by Peligrad, M. and Utev, S. \paper A new maximal inequality
and invariance principle for stationary sequences \jour Ann. Probab.
\vol 33 \pages 798-815 \yr 2005 \endref
\ref \key R \by Rosenblatt, M. \book Markov Processes: Structure and asymptotic behavior \publ Springer \publaddr Berlin \yr 1971 
\endref
\ref \key Vo 1 \by Voln\'y, D. \paper Approximating martingales and the
central limit theorem for strictly stationary processes
\jour Stochastic Processes and their Applications \vol 44 \pages 41-74
\yr 1993 \endref
\ref \key Vo 2 \by Voln\'y, D. \paper Martingale approximation of non adapted stochastic processes
with nonlinear growth of variance \paperinfo Dependence in Probability and Statistics
Series: Lecture Notes in Statistics, Vol. 187 Bertail, Patrice; Doukhan, Paul; Soulier, Philippe (Eds.) 
\yr 2006 \endref
\ref \key Wo \by Woodroofe, M. \paper A central limit theorem for functions of a Markov chain with applications to
shifts \jour Stoch. Proc. and their Appl. \vol 41 \pages 31-42 \yr 1992 \endref
\ref  \key Wu-Wo \by Wu, W.B. and Woodroofe, M. \paper Martingale approximation for 
sums of stationary processes \jour Ann. Probab. \vol 32 \pages 1674-1690 \yr 2004 \endref  
\endRefs
\enddocument
\end

\end